\documentclass[12pt,a4paper,twoside]{article}

\usepackage{amsmath, amssymb}
\usepackage{ascmac}

\DeclareMathOperator{\NS}{NS}
\DeclareMathOperator{\Pic}{Pic}

\DeclareMathOperator{\Cliff}{Cliff}

\DeclareMathOperator{\rk}{rk}
\begin{document}

\title{\textbf{\Large{A remark on the conjecture of Donagi-Morrison}}}

\author{Kenta Watanabe \thanks{Nihon University, College of Science and Technology,   7-24-1 Narashinodai Funabashi city Chiba 274-8501 Japan , {\it E-mail address:watanabe.kenta@nihon-u.ac.jp}}}

\date{}

\maketitle 

\begin{abstract}

\noindent Let $X$ be a K3 surface, let $C$ be a smooth curve of genus $g$ on $X$, and let $A$ be a base point free and primitive line bundle $g_d^r$ on $C$ with $d\geq4$ and $r\geq\sqrt{\frac{d}{2}}$. In this paper, we prove that if $g>2d-3+(r-1)^2$, then there exists a line bundle $N$ on $X$ which is adapted to $|C|$ such that $|A|$ is contained in the linear system $|N\otimes\mathcal{O}_C|$, and $\Cliff(N\otimes\mathcal{O}_C)\leq \Cliff(A)$.
 
\end{abstract}


\smallskip


\section{Introduction}

We work over the complex number field $\mathbb{C}$. Let $C$ be a smooth irreducible curve of genus $g$. Let $A$ be a line bundle $g_d^r$ on $C$. Then the following map $\mu_{0,A}$ defined by multiplication is called the {\it{Petri map}} associated with $A$.
$$\mu_{0,A}:H^0(A)\otimes H^0(K_C\otimes A^{\vee})\longrightarrow H^0(K_C).$$
Hence, $C$ is often said to be {\it{Brill-Noether-Petri}} ({\it{BNP}} for short), if the map $\mu_{0,A}$ is injective for any line bundle $A$ on $C$. It is well known that if $\mu_{0,A}$ is injective for $A$, then the Brill-Noether locus $W_d^r(C)$ consisting of line bundles of degree $d$ on $C$ which have at least $r+1$ linearly independent global sections is smooth at $A$. This means that if $C$ is BNP and $\rho(g,r,d)=g-(r+1)(g-d+r)>0$,  then $W_d^r(C)$ is smooth and irreducible away from $W_d^{r+1}(C)$, and has expected dimension $\rho(g,r,d)$. Here, $\rho(g,r,d)$ is called the {\it{Brill-Noether number}}. In particular, the following result is known.

\newtheorem{thm}{Theorem}[section]

\begin{thm} {\rm{([9, Theorem])}}. Let $C_0$ be a smooth connected curve lying on a K3 surface $X$. If any divisor in the linear system $|C_0|$ on $X$ is reduced and irreducible, then the general curve $C\in |C_0|$ is BNP. \end{thm}

From now on, let $X$ be a K3 surface, and assume that $C$ lies on $X$. Then $K_C\cong \mathcal{O}_C(C)$, and if $C$ is not hyperelliptic, then $C$ is a canonical curve (cf. [6], [13]). By Theorem 1.1, if the Picard group of $X$ is generated by the class of $C$, then the general curve which is linearly equivalent to $C$ is BNP. Hence, it is interesting to consider the problem of whether $C$ is BNP or not, in the case where the Picard number of $X$ is greater than one. Obviously, if $C$ is Brill-Noether special, then $C$ is not BNP. Hence, previously, the existence of a certain kind of lift of a line bundle on $C$ whose Brill-Noether number is negative has been investigated by several people (for example, see [1], [2], [4], [5], [10], [11], and [12]) to approach the above problem.

We set $L=\mathcal{O}_X(C)$. Then we say that a line bundle $N$ on $X$ is {\it{adapted}} to $|L|$ if the following conditions are satisfied.

\smallskip

\smallskip

(i) $h^0(N)\geq2$ and $h^0(L\otimes N^{\vee})\geq2$,

\smallskip

(ii) $h^0(N\otimes\mathcal{O}_{C^{'}})$ is independent of the smooth curve $C^{'}\in |L|$.

\smallskip

\smallskip

\noindent We can easily see that if a line bundle $N$ on $X$ satisfies the condition (i) as above, then $N\otimes\mathcal{O}_C$ contributes to the Clifford index of $C$ (see Notations and conventions for the definition of the Clifford index of a curve). However, the condition (ii) as above is not trivial due to the following exception.

$\;$

\noindent{\bf{Example 1.1}} ([12, Theorem 3.1, and Proposition 3.2]). There exists a smooth K3 surface $X\subset\mathbb{P}^5$ which contains a $(-2)$-curve $\Gamma$ and an elliptic curve $\Delta$ on $X$ with $\Delta.\Gamma=2$ such that all divisors of $|\Delta|$ are irreducible. Moreover, if we set $L=\mathcal{O}_X(a\Delta+(a-1)\Gamma)$ for an integer $a\geq3$, then $h^0(\mathcal{O}_C(2\Delta))$ depends on the smooth curve $C\in |L|$.

$\;$

\noindent Donagi and Morrison [3, Conjecture (1.2)] have conjectured the following.

\smallskip

\smallskip

\noindent{\bf{Conjecture 1.1}}. {\it{Let $X$ be a K3 surface, let $L$ be a base point free and big line bundle on $X$, and let $C\in |L|$ be a smooth curve of genus $g$. If a line bundle $A$ of degree $d\leq g-1$ on $C$ is base point free, and has negative Brill-Noether number, then there exists a line bundle $N$ on $X$ which satisfies the following conditions:

\smallskip

\smallskip

{\rm{(i)}} $|A|$ is contained in the restriction of $|N|$ to $C$;

\smallskip

{\rm{(ii)}} $\Cliff(N\otimes\mathcal{O}_C)\leq \Cliff(A)$}};

\smallskip

{\rm{(iii)}} $N.L\leq g-1$.

\smallskip

\smallskip

\noindent We often can not take a line bundle $N$ on $X$ satisfying the conditions (i) and (ii) as in Conjecture 1.1 to satisfy the condition (iii).

$\;$

\noindent{\bf{Example 1.2}} ([11, p554, Counterexample 1]). Let $X$ be a K3 surface, and let $B$ be a smooth curve of genus 2 on $X$. Assume that the Picard group of $X$ is generated by the class of $B$. Then there exist a smooth curve $C\in|3B|$, and a base point free line bundle $A=g_6^2$ on $C$ which is not isomorphic to $\mathcal{O}_C(B)$ such that $|A|$ is contained in the restriction of $|2B|$ to $C$. Obviously, $\mathcal{O}_X(2B)$ is adapted to $|C|$ and satisfies $\Cliff(\mathcal{O}_C(2B))=\Cliff(A)$. However, the genus of $C$ is 10 and $2B.C=12$.

$\;$

\noindent By the above reason, a line bundle $N$ on $X$ satisfying the conditions (i) and (ii) as in Conjecture 1.1 which is adapted to $|L|$ is often called a {\it{Donagi-Morrison lift}} of $A$. Previously, the above conjecture has been approached, under the hypothesis that $L$ is ample. For instance, Lelli-Chiesa [10, Theorem 1.1, and Corollary 1.3] and Rasmussen [12, Theorem 1.2] have investigated the existence of a Donagi-Morrison lift of $A$, in the case where $|A|$ is a net. On the other hand, if $C$ is neither hyperelliptic nor trigonal, and $A$ computes the Clifford index of $C$, then a necessary and sufficient condition for $A$ to have a Donagi-Morrison lift is known ([11, Theorem 1.1]).

However, if $L$ is not ample, the existence of a Donagi-Morrison lift of $A$ is still not known. In this paper, we investigate a sufficient condition for $A$ to have such a lift, in the case where $L$ is not necessarily ample. Our main result is the following.

\begin{thm} Let $X$ be a K3 surface, let $L$ be a base point free and big line bundle on $X$, and let $C\in |L|$ be a smooth curve of genus $g$. Let $A$ be a base point free and primitive line bundle $g_d^r$ on $C$ with $d\geq4$ and $r\geq\sqrt{\frac{d}{2}}$. If $g>2d-3+(r-1)^2$, then there exists a line bundle $N$ on $X$ which is adapted to $|L|$ such that $|A|\subset |N\otimes\mathcal{O}_C|$ and $\Cliff(N\otimes\mathcal{O}_C)\leq \Cliff(A)$.\end{thm}

\noindent{\bf{Remark 1.1}}. In Theorem 1.2, $|N\otimes\mathcal{O}_C|$ does not necessarily coincide with the restriction of $|N|$ to $C$. But, if $h^1(L\otimes N^{\vee})=0$, then the line bundle $N$ on $X$ satisfies the condition (i) as in Conjecture 1.1.

\smallskip

\smallskip

Our plan of this paper is as follows. In Section 2, we recall some results about line bundles on K3 surfaces. In Section 3, we recall several properties of the Lazarsfeld-Mukai bundle associated with a smooth curve on a K3 surface and a base point free line bundle on it. In Section 4, we prove Theorem 1.2. 

\smallskip

\smallskip

\noindent{\bf{Notations and conventions}}. In this paper, a curve and a surface are smooth and projective. Let $X$ be a curve or a surface. We denote the canonical line bundle of $X$ by $K_X$. We denote by $|L|$ the linear system associated with a divisor or a line bundle $L$ on $X$. We call a linear system of dimension two on $X$ a net. For a torsion free sheaf $E$ on $X$, we denote the rank of $E$, the dual of $E$, and the $i$-th Chern class of $E$, by $\rk(E)$, $E^{\vee}$, and $c_i(E)$, respectively.

A base point free line bundle $A$ on a curve $C$ is said to be {\it{primitive}} if $|K_C\otimes A^{\vee}|$ is also base point free. For a line bundle $A$ of degree $d$ on $C$, we denote the Clifford index of $A$ by $\Cliff(A):=d-2\dim|A|$. We say that a line bundle $A$ on a curve $C$ contributes to the Clifford index of $C$ if $h^0(A)\geq2$ and $h^1(A)\geq2$. Then the Clifford index of $C$ is the minimum value of the Clifford indices of such line bundles on $C$, and we denote it by $\Cliff(C)$. 

Let $X$ be a surface. Then we denote by $\Pic(X)$ and $\NS(X)$ the Picard group and the N${\rm{\acute{e}}}$ron-Severi group of $X$, respectively. We call the rank of $\NS(X)$ the Picard number of $X$. $X$ is called a regular surface if $h^1(\mathcal{O}_X)=0$. Moreover, if $K_X$ is trivial, then we call $X$ a K3 surface. Note that if $X$ is a regular surface, then $\Pic(X)$ and $\NS(X)$ are isomorphic. Hence, the Picard number of a K3 surface $X$ is the rank of $\Pic(X)$. In this paper, we will write $\Pic(X)_{\mathbb{R}}:=\Pic(X)\otimes_{\mathbb{Z}}\mathbb{R}$.

\section{Preliminaries}

In this section, we recall some fundamental facts concerning line bundles on K3 surfaces. Let $X$ be a K3 surface, and let $E$ be a vector bundle on $X$. First of all, by the Riemann-Roch theorem, we obtain the following equality.

$$\chi(E)=2\rk(E)+\frac{c_1(E)^2}{2}-c_2(E),$$
where $\chi(E)=h^0(E)-h^1(E)+h^2(E)$. Since $K_X$ is trivial, by the Serre duality, we have
$$h^i(E)=h^{2-i}(E^{\vee})\; (0\leq i\leq 2).$$

\smallskip

\smallskip

\noindent{\bf{Remark 2.1}}. If a line bundle $L$ on $X$ satisfies $L^2\geq-2$, then we have $\chi(L)\geq1$, and hence, exactly one of  $h^0(L)$ or $h^0(L^{\vee})$ is positive.

$\;$

\noindent {\bf{Remark 2.2}}. By the adjunction formula, for a non-zero effective divisor $D$ on $X$, we have $D^2=2P_a(D)-2$, where $P_a(D)$ is the arithmetic genus of $D$. In particular, if $D$ is a smooth rational curve on $X$, then $D^2=-2$.

\newtheorem{df}{Definition}[section]

\begin{df} We call an effective divisor $D$ on $X$ satisfying $D^2=-2$, a {\rm{$(-2)$-divisor}} on $X$. \end{df}

\noindent It is well known that if a big line bundle $L$ (i.e., a line bundle $L$ satisfying $L^2>0$ and $h^0(L)>0$) on $X$ satisfies $h^1(L)\neq0$, then there exists a smooth rational curve $\Gamma$ on $X$ such that $L.\Gamma\leq-1$, by the Kawamata-Viehweg vanishing theorem. However, a big line bundle $L$ on $X$ with $h^1(L)=0$ is not necessarily nef. Here, we recall the following result.

\newtheorem{prop}{Proposition}[section]

\smallskip

\smallskip

\begin{prop} {\rm{([8, Theorem])}}. Let $L$ be a big line bundle on $X$. Then  $h^1(L)\neq0$ if and only if  there exists a $(-2)$-divisor $\Gamma$ on $X$ with $L.\Gamma\leq-2$.\end{prop}

\smallskip

\smallskip

\noindent Next we recall the following characterization of base point free line bundles on $X$, due to Saint-Donat.

\smallskip

\smallskip

\begin{prop} {\rm{([13, Proposition 2.6])}}. Let $L$ be a non-trivial line bundle on $X$. If $h^0(L)>0$ and $|L|$ has no fixed component, then one of the following cases occurs. 

\smallskip

\smallskip

\noindent {\rm{(i)}} $L^2>0$ and the general member of $|L|$ is a smooth irreducible curve of genus $\frac{L^2}{2}+1$. In this case, $h^1(L)=0$.

\noindent {\rm{(ii)}} If $L^2=0$, then there exist an elliptic curve $\Delta$ on $X$ and an integer $r\geq1$ satisfying $L\cong\mathcal{O}_X(r\Delta)$. In this case, $h^0(L)=r+1$ and $h^1(L)=r-1$. \end{prop}

\smallskip

\smallskip

\noindent If $C$ is an irreducible curve on $X$ with $C^2\geq0$, then $|C|$ has no base point ([13, Theorem 3.1]). Proposition 2.2 implies the following proposition.

\smallskip

\smallskip

\begin{prop} {\rm{{([13, Corollary 3.2])}}}. Let $L$ be a non-trivial line bundle with $h^0(L)>0$ on $X$. Then $|L|$ has no base point outside of its fixed components.  \end{prop}

\smallskip

\smallskip

\noindent A non-zero effective divisor $D$ on $X$ is said to be {\it{numerically $m$-connected}} if $D_1.D_2\geq m$, for any decomposition $D=D_1+D_2$ with non-zero effective divisors $D_1$ and $D_2$ on $X$. 

$\;$

\noindent{\bf{Remark 2.3}}. If an effective divisor $D$ on $X$ is 1-connected, then, by the Serre duality and the exact sequence
$$0\longrightarrow\mathcal{O}_X(-D)\longrightarrow\mathcal{O}_X\longrightarrow\mathcal{O}_D\longrightarrow0,$$
we have $h^1(\mathcal{O}_X(D))=0$. Conversely, any member of the linear system $|L|$ associated with a non-trivial line bundle $L$ on $X$ satisfying $h^0(L)>0$ and $h^1(L)=0$ is 1-connected. However, the linear system associated with a 1-connected effective divisor on $X$ is not necessarily base point free.

$\;$

\noindent We recall the following proposition used in Section 4.

\smallskip

\smallskip

\begin{prop} {\rm{{([13, Lemma 3.7])}}}. Let $L$ be a base point free and big line bundle on $X$. Then any member of $|L|$ is 2-connected.\end{prop}

\smallskip

\smallskip

\noindent{\bf{Remark 2.4}}. If two base point free line bundles $L_1$ and $L_2$ satisfy $L_1.L_2>0$, by Proposition 2.4, we have $L_1.L_2\geq2$.

$\;$

\noindent Finally, we recall the following notion concerning the intersection theory of K3 surfaces to mention an important fact for the proof of our main theorem. 

\smallskip

\smallskip

\begin{df} {\rm{([6, Chapter 8, Definition 1.1])}}. The {\rm{positive cone}} 
$$\mathcal{C}_X\subset\Pic(X)_{\mathbb{R}}$$
is the connected component of the set $\{\alpha\in\Pic(X)_{\mathbb{R}}\;|\; \alpha^2>0\}$ which  contains one ample class. \end{df}

\smallskip

\smallskip

\noindent{\bf{Remark 2.5}}. The set $\{\alpha\in\Pic(X)_{\mathbb{R}}\;|\; \alpha^2>0\}$ is the disjoint union of $\mathcal{C}_X$ and $-\mathcal{C}_X$. Hence, any big line bundle on $X$ is contained in $\mathcal{C}_X$. Thus, any two big line bundles $L_1$ and $L_2$ on $X$ satisfy $L_1.L_2>0$ ([6, Chapter 8, p 146]).

\section{Lazarsfeld-Mukai bundles on K3 surfaces}

In this section, we recall the definition and several properties of Lazarsfeld-Mukai bundles used in Section 4. Let $X$ be a K3 surface, and let $C$ be a smooth curve of genus $g\geq2$ on $X$. Let $A$ be a base point free line bundle of degree $d$ on $C$ with $\dim|A|\geq2$, and let $V$ be a subspace of $ H^0(A)$ which forms a base point free linear system $(A,V)$ of dimension $r\geq1$ on $C$. Then we denote the dual of the kernel of the evaluation map $ev : V\otimes\mathcal{O}_X\longrightarrow A$, by $E_{C,(A,V)}$. Then it is often called the {\it{Lazarsfeld-Mukai}} ({\it{LM}} for short) bundle on $X$ associated with $C$ and $(A,V)$. Since $(A,V)$ has no base point, the evaluation map $ev$ is surjective, and hence, $E_{C,(A,V)}$ is locally free. By taking the dual of the exact sequence associated with the evaluation map $ev$, we obtain the following exact sequence.
$$0\longrightarrow V^{\vee}\otimes\mathcal{O}_X\longrightarrow E_{C,(A,V)}\longrightarrow K_C\otimes A^{\vee}\longrightarrow 0.$$
If $V=H^0(A)$, then we will write $E_{C,(A,V)}=E_{C,A}$.

\smallskip

\smallskip

\begin{prop} {\rm{([11, Proposition 2.1])}}. Any LM bundle $E_{C,(A,V)}$ on $X$ has the following properties:

\smallskip

\smallskip

\noindent {\rm{(i)}} $\rk E_{C,(A,V)}=r+1$, $c_1(E_{C,(A,V)})=\mathcal{O}_X(C)$, and $c_2(E_{C,(A,V)})=d$;

\noindent {\rm{(ii)}} $h^1(E_{C,(A,V)})=h^0(A)-r-1$, and $h^2(E_{C,(A,V)})=0$;

\noindent {\rm{(iii)}} $E_{C,(A,V)}$ is globally generated off the set of base points of $|K_C\otimes A^{\vee}|$;

\noindent {\rm{(iv)}} $\chi(E_{C,(A,V)}^{\vee}\otimes E_{C,(A,V)})=2(1-\rho(g,r,d))$. \end{prop}

\smallskip

\smallskip

\noindent{\bf{Remark 3.1}}. The assertion (iii) as in Proposition 3.1 means that if $A$ is primitive, then $E_{C,(A,V)}$ is generated by its global sections. In particular, if $A$ computes the Clifford index of $C$, then $K_C\otimes A^{\vee}$ also does so, and hence, $E_{C,(A,V)}$ is globally generated. Moreover, the assertion (ii) implies that if $V=H^0(A)$, then $h^i(E_{C,A})$ vanishes for $i=1,2$. Lelli-Chiesa [11, Definition 1] generalized the definition of a LM bundle on $X$ based on the properties (ii) and (iii) as in Proposition 3.1 as follows.

\smallskip

\smallskip

\begin{df} Let $E$ be a torsion free sheaf on $X$ with $h^2(E)=0$. If one of the following conditions is satisfied, $E$ is called a {\rm{generalized Lazarsfeld-Mukai}} {\rm{(g.LM}} for short{\rm{) }}bundle on $X$.

\smallskip

\smallskip

\noindent {\rm{(i)}} $E$ is a locally free sheaf which is globally generated off a finite set.

\noindent {\rm{(ii)}} $E$ is globally generated. \end{df}

\noindent In particular, a g.LM bundle $E$ on $X$ with $c_1(E)^2=0$ has the following properties.

\begin{prop} {\rm{([11, Proposition 2.7])}}. Let $E$ be a g.LM bundle on $X$ with $c_1(E)^2=0$. Then $E$ is a locally free sheaf with $c_2(E)=0$ generated by its global sections. Moreover, if $h^1(E)=0$, then there exists an elliptic curve $\Delta$ on $X$ satisfying $E=\mathcal{O}_X(\Delta)^{\oplus\rk E}$. \end{prop}

\noindent The assertion of (i) as in Proposition 3.1 means that $\Cliff(A)=c_2(E_{C,A})-2(\rk E_{C,A}-1)$. Thus, we can generalize the notion of the Clifford index of a line bundle on a smooth curve on $X$ as follows.

\begin{df} {\rm{([11, Definition 2])}}. We denote the Clifford index of a g.LM bundle $E$ on $X$ by $\Cliff(E):=c_2(E)-2(\rk E-1)$. \end{df}

\noindent Finally, we recall the following proposition used in Section 4 concerning the Clifford index of a g.LM  bundle $E$ on $X$ with $c_1(E)^2>0$.

\begin{prop} {\rm{([11, Proposition 2.4, and Corollary 2.5])}}. Let $E$ be a g.LM bundle on $X$ with $c_1(E)^2>0$. Then $\Cliff(E)\geq0$. \end{prop}

\section{Proof of Theorem 1.2}

Let $X$ be a K3 surface, let $L$ be a base point free and big line bundle on $X$, and let $C\in |L|$ be a smooth curve of genus $g$. Let $A$ be a base point free line bundle on $C$ of degree $d\geq 4$ such that $\dim|A|=r\geq\sqrt{\frac{d}{2}}$. Assume that $g>2d-3+(r-1)^2$. Then we note that $\rho(g,r,d)<0$ and $d<g-1$.

First of all, we prepare the following propositions to prove Theorem 1.2.

$\;$

\begin{prop} Let $\Delta$ be an elliptic curve on $X$, and let $\Gamma$ be an effective divisor on $X$ such that $\Delta.\Gamma=1$ and $C.\Gamma\leq r-2$. Then there exists a smooth rational curve $\Gamma_1\subset \Gamma$ such that $\Delta.\Gamma_1=1$ and $C.\Gamma_1\leq r-2$. \end{prop}

{\it{Proof}}. Since  $\Delta.\Gamma=1$, $\Gamma$ is not zero. If we take an irreducible component $\Gamma_1$ of $\Gamma$ with $\Delta.\Gamma_1=1$, by Remark 2.4, $\Gamma_1$ is a smooth rational curve on $X$. On the other hand, since $C.\Gamma\leq r-2$ and $C$ is nef, we have $C.\Gamma_1\leq r-2$. $\hfill\square$

\begin{prop} Assume that there exists an elliptic curve $\Delta$ on $X$ with $A=\mathcal{O}_C(r\Delta)$. Then there exists a line bundle $N$ on $X$ satisfying $A\subset N\otimes\mathcal{O}_C$ and $\Cliff(N\otimes\mathcal{O}_C)\leq \Cliff(A)$ such that $N$ is adapted to $|L|$.\end{prop}

\noindent Here, we prove the following three lemmas.

\newtheorem{lem}{Lemma}[section]

\begin{lem} Assume that there exist an elliptic curve $\Delta$ on $X$ satisfying $A=\mathcal{O}_C(r\Delta)$, and a smooth rational curve $\Gamma_1$ on $X$ with $\Delta.\Gamma_1=1$ and $C.\Gamma_1\leq r-2$. Then $h^0(\mathcal{O}_X(r\Delta+k\Gamma_1))>2$ and $h^0(L\otimes\mathcal{O}_X(-r\Delta-k\Gamma_1))>2$ for a non-negative integer $k\leq [\frac{r}{2}]$. Moreover, if there exists a $(-2)$-divisor $\Gamma_2$ on $X$ satisfying $\Delta.\Gamma_2=1$, $0\leq\Gamma_1.\Gamma_2$, and $C.\Gamma_2\leq r-2$, then $h^0(\mathcal{O}_X(r\Delta+\Gamma_1+\Gamma_2))>2$ and $h^0(L\otimes\mathcal{O}_X(-r\Delta-\Gamma_1-\Gamma_2))>2$. \end{lem}

{\it{Proof}}. By Proposition 2.2, we have 
$$h^0(\mathcal{O}_X(r\Delta+k\Gamma_1))\geq h^0(\mathcal{O}_X(r\Delta))=r+1\geq3.$$
Since $g>2d-3+(r-1)^2$, $\Delta.\Gamma_1=1$, and $C.\Gamma_1\leq r-2$, by Remark 2.2, we obtain
$$(C-r\Delta-k\Gamma_1)^2>2d-6+2(r-1)^2-2(k-1)^2,$$
and 
$$C.(C-r\Delta-k\Gamma_1)>3d-8+2(r-1)^2-k(r-2).$$
Since $d\geq4$, $r\geq2$, and $0\leq k\leq[\frac{r}{2}]$, the right-hand sides of the above inequalities are positive. Since $C$ is nef, by the Riemann-Roch theorem, we have $h^0(L\otimes\mathcal{O}_X(-r\Delta-k\Gamma_1))>2$.

On the other hand, by Proposition 2.2, we have
$$h^0(\mathcal{O}_X(r\Delta+\Gamma_1+\Gamma_2))\geq h^0(\mathcal{O}_X(r\Delta))=r+1\geq3.$$
Since $d<g-1$, $0\leq\Gamma_1.\Gamma_2$, and $(C-r\Delta).\Gamma_i\leq-2$ for $i=1,2$, by Remark 2.2, we have $(C-r\Delta-\Gamma_1-\Gamma_2)^2\geq 2g-2-2d+4>4.$
Since $g>2d-3+(r-1)^2$, $d\geq4$, $r\geq2$, and $C.\Gamma_i\leq r-2$ for $i=1,2$, we have 
$$C.(C-r\Delta-\Gamma_1-\Gamma_2)>3d-4+2(r-1)^2-2r\geq6.$$
By the same reason as above, we get $h^0(L\otimes\mathcal{O}_X(-r\Delta-\Gamma_1-\Gamma_2))>4$. $\hfill\square$

$\;$

\begin{lem} Assume that there exist an elliptic curve $\Delta$ on $X$ satisfying $A=\mathcal{O}_C(r\Delta)$, and a smooth rational curve $\Gamma_1$ on $X$ with $\Delta.\Gamma_1=1$ and $C.\Gamma_1=r-2$. Then there exists a line bundle $N$ on $X$ satisfying $A\subset N\otimes\mathcal{O}_C$ and $\Cliff(N\otimes\mathcal{O}_C)\leq \Cliff(A)$ such that $N$ is adapted to $|L|$.\end{lem}

{\it{Proof}}. First of all, we assume that $h^1(L\otimes\mathcal{O}_X(-r\Delta-\Gamma_1))=0$. Since $r\Delta+\Gamma_1$ is 1-connected, by Remark 2.3, we have $h^1(\mathcal{O}_X(r\Delta+\Gamma_1))=0$. By Remark 2.1 and Lemma 4.1, we have $h^0(\mathcal{O}_X(r\Delta+\Gamma_1)\otimes L^{\vee})=0$. Let $C^{'}$ be any smooth member of $|L|$. Since we have the exact sequence
$$0\longrightarrow\mathcal{O}_X(r\Delta+\Gamma_1)\otimes L^{\vee}\longrightarrow\mathcal{O}_X(r\Delta+\Gamma_1)\longrightarrow\mathcal{O}_{C^{'}}(r\Delta+\Gamma_1)\longrightarrow0,$$
we obtain $h^0(\mathcal{O}_{C^{'}}(r\Delta+\Gamma_1))=h^0(\mathcal{O}_X(r\Delta+\Gamma_1))=r+1$. On the other hand, by Remark 2.1 and Lemma 4.1, we have $h^0(\mathcal{O}_X(r\Delta)\otimes L^{\vee})=0$. By Proposition 2.2 and the exact sequence
$$0\longrightarrow\mathcal{O}_X(r\Delta)\otimes L^{\vee}\longrightarrow\mathcal{O}_X(r\Delta)\longrightarrow\mathcal{O}_{C^{'}}(r\Delta)\longrightarrow0,\leqno (1)$$
we have 
$$r+1=h^0(\mathcal{O}_X(r\Delta))\leq h^0(\mathcal{O}_{C^{'}}(r\Delta))\leq h^0(\mathcal{O}_{C^{'}}(r\Delta+\Gamma_1)).$$ Hence, we obtain $h^0(\mathcal{O}_{C^{'}}(r\Delta))=r+1$. By Lemma 4.1, $\mathcal{O}_X(r\Delta)$ is adapted to $|L|$, and $\Cliff(\mathcal{O}_C(r\Delta))=\Cliff(A)$. 

Assume that $h^1(L\otimes\mathcal{O}_X(-r\Delta-\Gamma_1))\neq 0$. Since $(C-r\Delta).\Gamma_1=-2$, $d\geq4$, $r\geq2$, and $g>2d-3+(r-1)^2$, we have
$$(C-r\Delta-\Gamma_1)^2=2g-2d>2d-6+2(r-1)^2>2(r-1)^2>0.\leqno (2)$$
Hence, there exists a $(-2)$-divisor $\Gamma_2$ on $X$ such that $(C-r\Delta-\Gamma_1).\Gamma_2\leq -2$, by Lemma 4.1 and Proposition 2.1.

We show that $\Delta.\Gamma_2=1$. We note that $\Delta.\Gamma_2\geq0$. Assume that $\Delta.\Gamma_2=0$. Since $(C-\Gamma_1).\Gamma_2\leq-2$, we have $\Gamma_1.\Gamma_2\geq C.\Gamma_2+2\geq2$. This means that $(\Delta+\Gamma_1+\Gamma_2)^2\geq2$. Hence, by Remark 2.5, both of the classes of $\Delta+\Gamma_1+\Gamma_2$ and $C-r\Delta-\Gamma_1$ in $\Pic(X)$ are contained in the positive cone, and we have
$$0<(C-r\Delta-\Gamma_1).(\Delta+\Gamma_1+\Gamma_2)=\frac{d}{r}-1+(C-r\Delta-\Gamma_1).\Gamma_2\leq\frac{d}{r}-3.$$
By the Hodge index theorem, we obtain $2(r-1)<\frac{d}{r}-3$. This contradicts the hypothesis that $r\geq\sqrt{\frac{d}{2}}$. Therefore, we have $\Delta.\Gamma_2\geq1$. If $\Delta.\Gamma_2\geq2$, then $(\Delta+\Gamma_2)^2\geq2$. By the same reason as above, we get
$$0<(C-r\Delta-\Gamma_1).(\Delta+\Gamma_2)=\frac{d}{r}-1+(C-r\Delta-\Gamma_1).\Gamma_2\leq\frac{d}{r}-3.$$
This gives the same contradiction as above, and hence, we obtain $\Delta.\Gamma_2=1$. 

We show that $\Gamma_1.\Gamma_2=0$. Assume that $\Gamma_1.\Gamma_2\leq -1$. Then $\Gamma_1$ is an irreducible component of $\Gamma_2$. Since $(C-r\Delta-\Gamma_1).\Gamma_2\leq-2$, we have $(C-r\Delta).\Gamma_2\leq\Gamma_1.\Gamma_2-2\leq-3$. Since $\Gamma_2-\Gamma_1$ is an effective divisor, $C.\Gamma_1\leq C.\Gamma_2\leq r-3$. This contradicts the hypothesis that $C.\Gamma_1=r-2$. Assume that $\Gamma_1.\Gamma_2\geq1$. Then we have $(\Delta+\Gamma_1+\Gamma_2)^2=2\Gamma_1.\Gamma_2\geq2$. Therefore, by the inequality (2) and the same reason as above, we have a contradiction. Hence, we have $\Gamma_1.\Gamma_2=0$. By the above argument, we have $C.\Gamma_2\leq r-2$. 

We consider the case where $h^0(\mathcal{O}_X(\Gamma_2-\Gamma_1))=0$. By Proposition 4.1, there exists a smooth rational curve $\Gamma_2^{'}\subset \Gamma_2$ with $\Delta.\Gamma_2^{'}=1$ and $C.\Gamma_2^{'}\leq r-2$. We set $N=\mathcal{O}_X(r\Delta+\Gamma_1+\Gamma_2^{'})$. Since $\Gamma_1\neq\Gamma_2^{'}$, we have $\Gamma_1.\Gamma_2^{'}\geq0$. Hence, by Lemma 4.1, we have $h^0(L\otimes N^{\vee})>2$. By Remark 2.1, we have $h^0(N\otimes L^{\vee})=0$. Since $r\Delta+\Gamma_1+\Gamma_2^{'}$ is 1-connected, $h^1(N)=0$. Hence, by the  Riemann-Roch theorem, we have $h^0(N)=\chi(N)\geq2r$. Since, for any smooth member $C^{'}\in |L|$, we have the exact sequence
$$0\longrightarrow N\otimes L^{\vee}\longrightarrow N\longrightarrow N\otimes\mathcal{O}_{C^{'}}\longrightarrow0,$$
we obtain $h^0(N\otimes\mathcal{O}_{C^{'}})=h^0(N)+h^1(N\otimes L^{\vee})$. Therefore, $N$ is adapted to $|L|$. Moreover, we have $h^0(N\otimes\mathcal{O}_C)\geq2r$. Since $r\Delta.C=d$, $C.\Gamma_1=r-2$, and $C.\Gamma_2^{'}\leq r-2$, we have $\Cliff(N\otimes\mathcal{O}_C)<\Cliff(A)$. 

We consider the case where $h^0(\mathcal{O}_X(\Gamma_2-\Gamma_1))>0$. Let $D$ be the movable part of $|\Gamma_2|$. Since $\Delta$ is nef, we have $\Delta.D\leq \Delta.\Gamma_2=1$. Assume that $D\neq\emptyset$. By Remark 2.4, we have $\Delta.D=0$. Moreover, by Proposition 2.4, we have $D^2=0$. Since $(D-\Delta)^2=0$, by Remark 2.1, we obtain $\mathcal{O}_X(\Delta)\subset\mathcal{O}_X(D)$. Since $\Delta.\Gamma_1=1$, $\Gamma_1$ is not an irreducible component of any member of $|\Delta|$. Hence, we have $\mathcal{O}_X(\Delta+\Gamma_1)\subset \mathcal{O}_X(\Gamma_2)$. Therefore, we have $\frac{d}{r}+r-2=C.(\Delta+\Gamma_1)\leq C.\Gamma_2\leq r-2$. This contradicts the fact that $\frac{d}{r}>0$. Thus $D=\emptyset$.

Here, let $\Gamma_2^{'}\subset \Gamma_2$ be the minimal effective divisor on $X$ with $\Gamma_1.\Gamma_2^{'}=0$, $\Delta.\Gamma_2^{'}=1$, and ${\Gamma_2^{'}}^2=-2$. We set $N=\mathcal{O}_X(r\Delta+\Gamma_1+\Gamma_2^{'})$. Then, by Lemma 4.1, we have $h^0(L\otimes N^{\vee})>2$. We show that $N$ is nef. Assume that there exists a smooth rational curve $\Gamma$ on $X$ with $N.\Gamma<0$. Since $N.\Gamma_1=r-2\geq0$, we have $\Gamma\neq\Gamma_1$. Hence, we have $\Gamma_1.\Gamma\geq0$. Since $\Delta$ is nef, we have $\Gamma_2^{'}.\Gamma\leq N.\Gamma<0$. This means that $\Gamma$ is an irreducible component of $\Gamma_2^{'}$. Note that since $N.\Gamma_2^{'}=r-2\geq0$, we have $\Gamma_2^{'}\neq \Gamma$.

First of all, we obtain $\Gamma_1.\Gamma=0$ and $\Gamma_2^{'}.\Gamma=-1$. In fact, since 
$$\Gamma_2^{'}.\Gamma=(N-\Gamma_1-r\Delta).\Gamma\leq N.\Gamma-\Gamma_1.\Gamma,$$
if $\Gamma_1.\Gamma\geq1$ or $\Gamma_2^{'}.\Gamma\leq -2$, we have $(\Gamma_2^{'}-\Gamma).\Gamma\leq0$. Since $\Gamma_2^{'}-\Gamma$ is a non-zero effective divisor on $X$, this means that $\Gamma_2^{'}$ is not 1-connected. Since $h^0(\mathcal{O}_X(\Gamma_2))=1$, we have $h^0(\mathcal{O}_X(\Gamma_2^{'}))=1$. By the Riemann-Roch theorem, we have $h^1(\mathcal{O}_X(\Gamma_2^{'}))=0$. By Remark 2.3, this is a contradiction. 

On the other hand, since $\Delta.\Gamma_1=\Delta.\Gamma_2=1$, we have $\Delta.(\Gamma_2-\Gamma_1)=0$. Since $h^0(\mathcal{O}_X(\Gamma_2-\Gamma_1))>0$ and $D=\emptyset$, $\Gamma_1$ is an irreducible component of $\Gamma_2$. Since $\Gamma\subset\Gamma_2-\Gamma_1$, we obtain $\Delta.\Gamma=0$.

By the above argument, we obtain $\Gamma_1.(\Gamma_2^{'}-\Gamma)=0$, $\Delta.(\Gamma_2^{'}-\Gamma)=1$, and $(\Gamma_2^{'}-\Gamma)^2=-2$. This contradicts the minimality of $\Gamma_2^{'}$. Therefore, $N$ is nef. Since $N^2=4r-4>0$, by the Kawamata-Viehweg vanishing theorem, we have $h^1(N)=0$. Hence, by the Riemann-Roch theorem, we have $h^0(N)=\chi(N)=2r$. Therefore, $N$ is adapted to $|L|$. Since $h^0(N\otimes\mathcal{O}_C)\geq2r$. By the same reason as above, we obtain $\Cliff(N\otimes\mathcal{O}_C)<\Cliff(A)$. $\hfill\square$

\begin{lem} Assume that there exist an elliptic curve $\Delta$ on $X$ satisfying $A=\mathcal{O}_C(r\Delta)$, and a smooth rational curve $\Gamma_1$ on $X$ with $\Delta.\Gamma_1=1$ and $C.\Gamma_1\leq r-3$. Then there exists a line bundle $N$ on $X$ satisfying $A\subset N\otimes\mathcal{O}_C$ and $\Cliff(N\otimes\mathcal{O}_C)\leq \Cliff(A)$ such that $N$ is adapted to $|L|$.\end{lem}

{\it{Proof}}. Since $C.\Gamma_1\geq0$, we have $r\geq3$. Assume that $r=3$. Then we obtain $C.\Gamma_1=0$. By Remark 2.1 and Lemma 4.1, $h^0(\mathcal{O}_X(3\Delta+\Gamma_1)\otimes L^{\vee})=0$. Since $3\Delta+\Gamma_1$ is 1-connected, by Remark 2.3, we have $h^1(\mathcal{O}_X(3\Delta+\Gamma_1))=0$. Hence, by the Riemann-Roch theorem, $h^0(\mathcal{O}_X(3\Delta+\Gamma_1))=4$. Since $h^0(A)=4$, by the exact sequence
$$0\longrightarrow\mathcal{O}_X(3\Delta+\Gamma_1)\otimes L^{\vee}\longrightarrow\mathcal{O}_X(3\Delta+\Gamma_1)\longrightarrow A\longrightarrow0,$$
we have $h^1(\mathcal{O}_X(3\Delta+\Gamma_1)\otimes L^{\vee})=0$. For any smooth curve $C^{'}\in |L|$, we obtain the exact sequence
$$0\longrightarrow\mathcal{O}_X(3\Delta+\Gamma_1)\otimes L^{\vee}\longrightarrow\mathcal{O}_X(3\Delta+\Gamma_1)\longrightarrow \mathcal{O}_{C^{'}}(3\Delta)\longrightarrow0.$$
Hence, we have $h^0(\mathcal{O}_{C^{'}}(3\Delta))=4$. This implies that $\mathcal{O}_X(3\Delta)$ is adapted to $|L|$, and $\Cliff(\mathcal{O}_C(3\Delta))=\Cliff(A)$. 

Assume that $r\geq4$. We set $N=\mathcal{O}_X(r\Delta+2\Gamma_1)$. By Lemma 4.1, we have $h^0(L\otimes N^{\vee})>2$.  Since $N$ is nef and big, by the Kawamata-Viehweg vanishing theorem, we have $h^1(N)=0$. By the Riemann-Roch theorem, $h^0(N)=2r-2\geq6$. Hence, $N$ is adapted to $|L|$. Since $C.\Gamma_1\leq r-3$, we have $\Cliff(N\otimes\mathcal{O}_C)\leq\Cliff(A)$. $\hfill\square$

$\;$

{\it{Proof of Proposition 4.2}}. Since $h^0(\mathcal{O}_X(r\Delta)\otimes L^{\vee})=0$, if $h^1(L\otimes\mathcal{O}_X(-r\Delta))=0$, by Proposition 2.2 and the exact sequence (1), we have $h^0(\mathcal{O}_{C^{'}}(r\Delta))=h^0(\mathcal{O}_X(r\Delta))=r+1$, for any smooth irreducible curve $C^{'}\in |L|$. In this case, by Lemma 4.1, $\mathcal{O}_X(r\Delta)$ is adapted to $|L|$, and $\Cliff(\mathcal{O}_C(r\Delta))=\Cliff(A)$. 

We assume that $h^1(L\otimes\mathcal{O}_X(-r\Delta))\neq0$. Since $g>2d-3+(r-1)^2$, $d\geq4$, and $r\geq2$, we have
$$(C-r\Delta)^2=2g-2-2d>2d-8+2(r-1)^2\geq 2(r-1)^2>0.$$
Hence, by Lemma 4.1 and Proposition 2.1, there exists a $(-2)$-divisor $\Gamma$ on $X$ such that $(C-r\Delta).\Gamma\leq -2$. Since $L$ and $\Delta$ are nef, we have $\Delta.\Gamma\geq1$. We show that $\Delta.\Gamma=1$. Assume that $\Delta.\Gamma\geq2$. Then we have $(\Delta+\Gamma)^2\geq2$. Hence, by Remark 2.5, both of the classes of $\Delta+\Gamma$ and $C-r\Delta$ in $\Pic(X)$ are contained in the positive cone, and we have
$$0<(C-r\Delta).(\Delta+\Gamma)=\frac{d}{r}+(C-r\Delta).\Gamma\leq\frac{d}{r}-2.$$
By the Hodge index theorem, we have $2(r-1)<\frac{d}{r}-2$. This contradicts the hypothesis that $r\geq\sqrt{\frac{d}{2}}$. Hence, $\Delta.\Gamma=1$ and $C.\Gamma\leq r-2$. By Proposition 4.1, Lemma 4.2, and Lemma 4.3, the assertion of Proposition 4.2 holds. $\hfill\square$

\smallskip

\smallskip

From now on, we assume that $A$ is primitive. Then the Lazarsfeld-Mukai bundle $E_{C,A}$ associated with $C$ and $A$ is globally generated by Remark 3.1. Although if $E_{C,A}$ contains a base point free and big line bundle $M$ on $X$ as a saturated subsheaf of $E_{C,A}$, then the consequence of Theorem 1.2 holds for $N:=L\otimes M^{\vee}$ due to Proposition 5.1 as in [11], we can construct such a sub-line bundle $M\subset E_{C,A}$ with $h^0(M)\geq2$ which is not necessarily base point free. Here, we obtain the following lemma by the same way of the proof of Lemma 4.1 as in [14]. 

\begin{lem} Let $E_{C,A}$ be as above. Then there exist a torsion free sheaf $F$ of rank $r$ on $X$ generated by its global sections, and a saturated line bundle $M\subset E_{C,A}$ such that $h^0(M)\geq2$ and $E_{C,A}$ sits in the following exact sequence.
$$0\longrightarrow M\longrightarrow E_{C,A}\longrightarrow F \longrightarrow0.\leqno (3)$$
\end{lem}

{\it{Proof of Theorem 1.2}}. The torsion free sheaf $F$ on $X$ as in Lemma 4.4 sits in the following exact sequence.
$$0\longrightarrow F\longrightarrow F^{\vee\vee}\longrightarrow S \longrightarrow0,\leqno (4)$$
where $S$ is a coherent sheaf of finite length on $X$. We set $c_1(F^{\vee\vee})=N$ and let $W$ be the support of $S$. Then, by [7, (0.3)], we have
$$c_2(E_{C,A})=\deg A=d=M.N+\ell(W)+c_2(F^{\vee\vee}),\leqno (5)$$
where $\ell(W)$ is the length of $W$. Since $h^2(E_{C,A})=0$, by the exact sequence (3), we have $h^2(F)=0$. Hence, by the exact sequence (4), we obtain $h^2(F^{\vee\vee})=0$. On the other hand, since $h^1(E_{C,A})=0$, by the exact sequence (3), we have $h^1(F)=0$. Since $F$ is globally generated,  $c_1(F^{\vee\vee})$ is base point free off the set $W$. By Proposition 2.3, $N$ is base point free everywhere.

On the other hand, by the exact sequence (3), we have $h^0(E_{C,A}\otimes M^{\vee})>0$. If we apply the functor $\otimes M^{\vee}$ to the exact sequence
$$0\longrightarrow H^0(A)^{\vee}\otimes\mathcal{O}_X\longrightarrow E_{C,A}\longrightarrow K_C\otimes A^{\vee}\longrightarrow0,$$
we obtain the exact sequence
$$0\longrightarrow H^0(A)^{\vee}\otimes M^{\vee}\longrightarrow E_{C,A}\otimes M^{\vee}\longrightarrow N\otimes\mathcal{O}_C\otimes A^{\vee}\longrightarrow0.$$
Since $h^0(M^{\vee})=0$, we have $h^0(N\otimes\mathcal{O}_C\otimes A^{\vee})\geq h^0(E_{C,A}\otimes M^{\vee})>0$. This means that 
$$A\subset N\otimes\mathcal{O}_C.\leqno (6)$$

Assume that $N^2=0$. Since $F$ is a g.LM bundle with $h^1(F)=0$, by Proposition 3.2, we can take an elliptic curve $\Delta$ on $X$ with $F=F^{\vee\vee}\cong\mathcal{O}_X(\Delta)^{\oplus r}$, and obtain $c_2(F^{\vee\vee})=0$ and $\ell(W)=0$. Hence, we have $L.N=M.N=d$. Since $N\cong\mathcal{O}_X(r\Delta)$, by (6), it turns out that $A\cong \mathcal{O}_C(r\Delta)$. By Proposition 4.2, we obtain the consequence of Theorem 1.2.

Assume that $N^2>0$. By Proposition 2.2, we have $h^1(N)=0$. Since $h^0(M^{\vee})=0$, and, for any smooth curve $C^{'}\in |L|$, we have the exact sequence
$$0\longrightarrow M^{\vee}\longrightarrow N\longrightarrow N\otimes\mathcal{O}_{C{'}}\longrightarrow0,$$
we obtain $h^0(N\otimes\mathcal{O}_{C{'}})=h^0(N)+h^1(M^{\vee})$. This means that $N$ is adapted to $|L|$. Moreover, since $F^{\vee\vee}$ is a g.LM bundle of type (i) as in Definition 3.1, we have $\Cliff(F^{\vee\vee})\geq0$ by Proposition 3.3. This implies that $c_2(F^{\vee\vee})\geq 2(r-1)$. Since $\ell(W)\geq0$, by the equality (5), we have

\smallskip

\smallskip

\noindent $\Cliff(N\otimes\mathcal{O}_C)=N.L-2(h^0(N\otimes\mathcal{O}_C)-1)$

$\leq N.L-2(h^0(N)-1)$

$=N.M-2=d-\ell(W)-c_2(F^{\vee\vee})-2$

$\leq d-2(r-1)-2=d-2r=\Cliff(A).$ $\hfill\square$

$\;$

\noindent{\bf{Example 4.1}} ([14, Proposition 5.4]). Let $C_0$ be a smooth plane curve of degree 4, and let $\pi:C\longrightarrow C_0$ be a double covering branched at distinct $24$ points on $C_0$ such that $C$ is contained in a K3 surface $X$. If $\pi$ is the morphism associated with a base point free net $A$ of degree 8 on $C$, then there exists a base point free line bundle $N$ on $X$ which is adapted to $|C|$ such that $N^2=2$, $A=N\otimes\mathcal{O}_C$, and $\Cliff(N\otimes\mathcal{O}_C)=\Cliff(A)$.

\smallskip

\smallskip

\noindent If $A$ computes the Clifford index of $C$, then $A$ is primitive. However, such a condition often can not be expected. Here, we recall the following notion.

\begin{df} {\rm{([2, Definition 1.16])}}. We call the minimum value of the Clifford indices of line bundles $A$ on a curve $C$ with $h^0(A)\geq2$, $h^1(A)\geq2$, and $\rho(A)<0$, the {\rm{special Clifford index}} of $C$, where $\rho(A)$ is the Brill-Noether number of $A$.\end{df}

\noindent We can easily see that if a line bundle $A$ on a curve $C$ computes the special Clifford index of $C$, then $K_C\otimes A^{\vee}$ also does so, and $A$ is primitive. Hence, we obtain the following.

$\;$

\noindent{\bf{Corollary 4.1}} Let $X$ be a K3 surface, let $L$ be a base point free and big line bundle on $X$, and let $C\in |L|$ be a smooth curve of genus $g$. Let $A$ be a line bundle $g_d^r$ on $C$ with $d\geq4$ and $r\geq\sqrt{\frac{d}{2}}$ which computes the special Clifford index of $C$. If $g>2d-3+(r-1)^2$, then there exists a line bundle $N$ on $X$ which is adapted to $|L|$ such that $|A|\subset |N\otimes\mathcal{O}_C|$ and $\Cliff(N\otimes\mathcal{O}_C)\leq \Cliff(A)$.

$\;$

%\noindent {\bf{Acknowledgements}}

\smallskip

\smallskip

%\noindent The author would like to thank the referee for very careful comments.

\end{document}